\documentclass[11pt,leqno]{article} 

\usepackage{fullpage}

\usepackage{amssymb} 
\usepackage{amsmath} 
\date{}
\title{Classification of linearly compact simple N=6 3-algebras} 
\author{{\sc Nicoletta Cantarini}\thanks{Dipartimento di Matematica
Pura ed Applicata, Universit\`a di Padova, Padova, Italy}
\and\setcounter{footnote}{6}
{\sc Victor G.\ Kac}\thanks{Department of Mathematics, MIT, Cambridge,
Massachusetts 02139, USA}}

\newtheorem{theorem}{Theorem}[section] 
\newtheorem{lemma}[theorem]{Lemma} 
 
\newtheorem{proposition}[theorem]{Proposition} 
\newtheorem{definition}[theorem]{Definition} 
\newtheorem{remark}[theorem]{Remark}

\newtheorem{example}[theorem]{Example}

\def\Z{\mathbb{Z}}
\def\ZZ{\mathbb{Z}}
\def\g{\mathfrak{g}}

\def\fg{\mathfrak{g}}

\def\C{\mathbb{C}}
\def\F{\mathbb{F}}
\def\FF{\mathbb{F}}

\def\ZZ{\mathbb{Z}}

\def\0{\bar{0}}
\def\1{\bar{1}}

\numberwithin{equation}{section}

\makeatletter

\def\enumerate{%
  \ifnum \@enumdepth >\thr@@\@toodeep\else
    \advance\@enumdepth\@ne
    \edef\@enumctr{enum\romannumeral\the\@enumdepth}%
      \list
        {\csname label\@enumctr\endcsname}%
        {\usecounter\@enumctr
          \addtolength{\leftmargin}{-\leftmargin}
          \settowidth{\labelwidth}{(99)}
          \itemindent = \labelwidth
          \addtolength{\itemindent}{\labelsep}
        \listparindent=1em      
          \def\makelabel##1{{##1}\hfill}
          }%
  \fi}

\makeatother

\begin{document} 
\maketitle 
\begin{abstract}
$N\leq 8$ $3$-algebras have recently appeared in $N$-supersymmetric 
$3$-dimensional Chern-Simons gauge theories. In our previous paper we
classified linearly compact simple $N=8$ $n$-algebras for any $n\geq 3$.
In the present paper we classify linearly compact simple $N=6$ $3$-algebras,
using their correspondence with simple linearly compact Lie superalgebras
with a consistent short $\ZZ$-grading, endowed with a graded conjugation.
We also briefly discuss $N=5$ $3$-algebras.
\end{abstract} 

\setcounter{section}{-1}
\section{Introduction}
\label{sec:intro}

In recent papers on N-supersymmetric 3-dimensional Chern-Simons
gauge theories various types of 3-algebras have naturally appeared (see
\cite{G}, \cite{BL1}, \cite{BL2}, \cite{ABJM}, \cite{BB}, $\ldots$).

Recall that a 3-algebra (also called a triple system \cite{J}) is a
vector space $V$ with a ternary (or 3-)bracket 
$V^{\otimes 3} \to V$, $a \otimes b \otimes c \mapsto [a,b,c]$.

The 3-algebras that appear in supersymmetric 3-dimensional
Chern-Simons theories satisfy certain symmetry conditions and a
Jacobi-like identity, usually called the fundamental identity
(very much like the Lie algebra bracket).

The simplest among them are 3-Lie algebras, for which  the
symmetry condition is the total anti-commutativity:
\begin{equation}
\label{eq:0.1}
   [a,b,c] =- [b,a,c] =- [a,c,b]\, , 
\end{equation}
and the fundamental identity is:
\begin{equation}
\label{eq:0.2}
    [a,b, [x,y,z]]  = [[a,b,x],y,z] + [x,[a,b,y],z]+[x,y,[a,b,z]]\, .
\end{equation}
(Of course, identity (0.2) simply says that for each $a,b \in V$,
the endomorphism $D_{a,b}(x) = [a,b,x]$ is a derivation of the
3-algebra $V$, very much like the Jacobi identity for Lie
algebras; in fact, this identity appears already in \cite{J}.)

The notion of a 3-Lie algebra generalizes to that of an $n$-Lie
algebra for an arbitrary integer $n \geq 2$ in the obvious way.
In this form they were introduced by Filippov in 1985 \cite{F}.  It
was subsequently proved in Ling's thesis \cite{L} that for each $n \geq 3$
there is only one simple finite-dimensional $n$-Lie algebra over
an algebraically closed field of characteristic~$0$,
by an analysis of the linear Lie algebra spanned by the derivations
$D_{a,b}$ (for $n=3$ this fact was independently proved in \cite{Fr}).  
This unique simple $n$-Lie algebra is the vector
product $n$-algebra $O^n$ in an $n+1$-dimensional vector space \cite{F},
\cite{L}.  Recall that, endowing an $n+1$-dimensional vector space $V$
with a non-degenerate symmetric bilinear form $(\cdot ,\cdot
)$ and choosing  a basis $\{ a_i \}$ and the dual basis $\{ a^i
\}$, so that $(a_i,a^j) = \delta_{ij}$, the vector product of
$n$~vectors from the basis $\{ a_i \}$ is defined by:
\begin{equation}
  \label{eq:0.3}
  [a_{i_1}, \ldots , a_{i_n}] = \epsilon_{i,\ldots, i_{n+1}}
     a^{i_{n+1}}\, ,
\end{equation}
where $\epsilon_{i_1,\ldots, i_{n+1}}$ is a non-zero totally
antisymmetric tensor, and extended by $n$-linearity.

In \cite{CantaK5} we extended this classification to linearly compact $n$-Lie
algebras.  Recall that a linearly compact $n$-algebra is a
topological $n$-algebra (i.e.,~the $n$-product is continuous) whose
underlying topological vector space is linearly compact.  The
basic examples of linearly compact spaces over a field $\F$ are the spaces of
formal power series $\FF [[x_1 , \ldots ,x_n]]$, endowed with
formal topology, or a direct sum of a finite number of such
spaces (these include finite-dimensional vector spaces with
discrete topology).  Our result is that a complete list of simple
linearly compact $n$-Lie algebras over an
algebraically closed field $\FF$ of characteristic $0$ for $n >2$
consists of four $n$-algebras:  the $n+1$-dimensional $O^n$ and
three infinite-dimensional, which we denoted by $S^n$, $W^n$,
$SW^n$ (see \cite{CantaK5} for their construction).

Our method consists in associating  to an 
$n$-Lie algebra $\fg$ a Lie superalgebra $L=
\bigoplus^{n-1}_{j=-1} L_j$ with a consistent $\ZZ$-grading,
such that $L_{-1} =\Pi \fg$  (i.e., the vector space ~$\fg$ with odd parity),
satisfying the following properties:
\begin{eqnarray}
  \label{eq:0.4}
&  [a,L_{-1}]=0\,,\, a \in L_j\,,\, j \geq 0\,,\,
 \hbox{\,\, imply \,\,} a=0\,  &\hbox{(transitivity)\,\,};\\
  \label{eq:0.5}
 &\dim L_{n-1}=1 \hbox{\,\, unless all $n$-brackets are\,\, } 0; &\\
%
  \label{eq:0.6}
 &[L_{j}, L_{-1}] = L_{j-1} \hbox{\,\, for all \,\,} j \, ; &\\
  \label{eq:0.7}
 & [L_j, L_{n-j-1}] = 0  \hbox{\,\, for all \,\,} j \, . &
\end{eqnarray}
Provided that $n>2$, the Lie superalgebra $[L,L]=
\bigoplus^{n-2}_{j=-1} L_j$ is simple, hence the classification
of simple linearly compact $n$-Lie algebras is thereby reduced to
the known classification of simple linearly compact Lie
superalgebras \cite{K1}, \cite{K3} and their consistent $\ZZ$-gradings
\cite{K2}, \cite{CantaK}.

Note that, given a consistently  $\ZZ$-graded Lie superalgebra
$L=\bigoplus^{n-1}_{j=-1} L_j$, satisfying
(\ref{eq:0.4})--(\ref{eq:0.7}), we can recover the $n$-bracket on
$\fg = \Pi L_{-1}$ by choosing a non-zero $\mu \in L_{n-1}$ and
letting:
\begin{displaymath}
  [a_1,\ldots ,a_n] = [\ldots [\mu ,a_1],\ldots , a_n]\, .
\end{displaymath}

The 3-Lie algebras appear in $N=8$ supersymmetric 3-dimensional
Chern-Simons theories, hence it is natural to call them the $N=8$
3-algebras.  The next in the hierarchy of 3-algebras are those
which appear in $N=6$ supersymmetric 3-dimensional Chern-Simons
theories (case $N=7$ reduces to $N=8$), which we shall call $N=6$
3-algebras.  They are defined by the following axioms:
\begin{eqnarray*}
  \label{eq:0.8}
     [a,b,c] &=&-[c,b,a]\, ;\\
  \label{eq:0.9}
     [a,b,[x,y,z]] &=& [[a,b,x],y,z]-[x,[b,a,y],z]+ [x,y,[a,b,z]]\, .
\end{eqnarray*}
Note that any $N=8$ 3-algebra is also an $N=6$ 3-algebra.  (It is
unclear how to define $N=6$ $n$-algebras for $n>3$.) 

The main goal of the present paper is the classification of simple
linearly compact $N=6$ 3-algebras over $\C$.  The method again consists of
associating to an $N=6$ 3-algebra a $\ZZ$-graded Lie
superalgebra, but in a different way (our construction in \cite{CantaK5}
uses total anti-commutativity in an essential way).

Given an $N=6$ 3-algebra $\fg$, following Palmkvist [P], we
associate to $\fg$ a pair $(L=L_{-1} \oplus L_0 \oplus L_1\, ,
\varphi)$, where $L$ is a consistently $\ZZ$-graded Lie
superalgebra with $L_{-1}=\Pi \fg$ and $\varphi$ is its
automorphism, such that the following properties hold:
\begin{eqnarray}
  \label{eq:0.10}
    && \hbox{transitivity \,\, (\ref{eq:0.4})} \, ;\\
  \label{eq:0.11}
  &&  [L_{-1},L_1] = L_0 \, ;\\
  \label{eq:0.12}
  &&   \varphi (L_j) = L_{-j} \hbox{\,\, and \,\,} \varphi^2 (a)=(-1)^j
        a \hbox{\,\, if \,\,} a \in L_j \, .
\end{eqnarray}
An automorphism $\varphi$ of the $\ZZ$-graded Lie superalgebra $L$,
satisfying (\ref{eq:0.12}), is called a {\em graded conjugation}.

It is easy to see that $\fg$ is a simple $N=6$ $3$-algebra if and only if
the associated Lie super algebra $L$ is simple.  Thus, this
construction reduces the classification of simple linearly
compact $N=6$ 3-algebras to the  classification of
consistent $\ZZ$-gradings of simple linearly compact Lie
superalgebras of the form $L=L_{-1} \oplus L_0 \oplus L_1$ (then
(\ref{eq:0.10}) and (\ref{eq:0.11}) automatically hold), and
their graded conjugations $\varphi$.  The 3-bracket on $\fg= \Pi L_{-1}$
is recovered by letting $[a,b,c]=[[a,\varphi (b)],c]$. 

The resulting classification of $N=6$ $3$-algebras is given by Theorem 1.11,
and its proof is given in Sections 2,3, and 4. In Section 5, based on the 
discussion in \cite {BB}, we propose a definition and obtain a classification
of $N=5$ $3$-algebras,
and give a construction of a reduction from $N=6$ $3$-algebras to $N=5$ 
$3$-algebras.

The base field is 
an algebraically closed field $\F$ of characteristic zero.

\section{Examples of $\mathbf{N=6}$ 3-algebras}
\begin{definition}\label{N6} An $N=6$ $3$-algebra is a $3$-algebra whose 
$3$-bracket $[\cdot,\cdot,\cdot]$ satisfies
the following axioms:
\begin{enumerate}
\item[$(a)$] $[u,v,w]=-[w,v,u]$
        \item[$(b)$]  $[u,v,[x,y,z]]=[[u,v,x],y,z]-[x,[v,u,y],z]+[x,y,[u,v,z]]$
\end{enumerate}
\end{definition}
\begin{example}\em Every 3-Lie algebra is an $N=6$ 3-algebra.
\end{example}


\begin{example}\label{sl(2,2)}\em
Let $A$ be an associative algebra and let ${}^*$ be an anti-involution of $A$, 
i.e., for every $a\in A$, $(a^*)^*=a$ and
for every $a, b\in A$, $(ab)^*=b^*a^*$.  
Then $A$ with $3$-bracket
\begin{equation}
[a,b,c]=ab^*c-cb^*a
\label{antitranspose}
\end{equation}
is an $N=6$ $3$-algebra.
For instance, if $A=M_n(\F)$ and ${}^*$ is the transposition map ${}^t: 
a\mapsto a^t$, then
the corresponding $3$-bracket (\ref{antitranspose}) defines on $A$ an 
$N=6$ 3-algebra structure.
Likewise, if $n=2k$ and ${}^*$ is the symplectic involution
$${}^*: a\mapsto J_{2k}a^tJ_{2k}^{-1},$$ 
with $J_{2k}=\left(
\begin{array}{cc}
0 & I_k\\
-I_k & 0
\end{array}
\right)$,
then (\ref{antitranspose}) defines on $M_{n}(\F)$ an $N=6$ 3-algebra structure.
If $n=2$, this is in fact a 3-Lie algebra structure.

More generally, consider the vector space $M_{m,n}(\F)$ of $m\times n$ 
matrices with entries in $\F$,
and let ${}^*: M_{m,n}(\F) \rightarrow M_{n,m}(\F)$ be a map satisfying the 
following property:
for $v,u,y\in M_{m,n}(\F)$,
\begin{equation}\label{*condition}
(vu^*y)^*=y^*uv^*.
\end{equation} 
Then $M_{m,n}(\F)$ with the $3$-bracket (\ref{antitranspose}) is an $N=6$ 
3-algebra.
(Note that if $m=n$, then property (\ref{*condition}) implies 
that ${}^*$ is an anti-involution).
As an example we can take ${}^*={}^t$; we denote this 3-algebra by 
$A^3(m,n;t)$.
Likewise, if $m=2h$ and $n=2k$, it is easy to check that
 the map $st: M_{m,n}(\F) \rightarrow M_{n,m}(\F)$, given by
\begin{equation}\label{symplectic}
 st: a\mapsto J_{2k}a^tJ_{2h}^{-1},
\end{equation} 
satisfies property (\ref{*condition}),
hence $M_{m,n}(\F)$ with the corresponding $3$-bracket (\ref{antitranspose}) is 
an $N=6$ 3-algebra, which we denote by $A^3(m,n;st)$. It is easy to see that
the 3-algebra $A^3(2,2;st)$ is actually an N=8 3-algebra, isomorphic to
$O^3$.
\end{example}

\begin{lemma}\label{sl(m,n)}
\begin{itemize}
\item[(a)] Let ${}^*: M_{m,n}(\F) \rightarrow M_{n,m}(\F)$ 
be defined by: $b^{*}=k^{-1}b^th$ for some symmetric matrices
$h\in GL_m(\F)$ and $k\in GL_n(\F)$, then the corresponding $3$-bracket 
(\ref{antitranspose}) is
isomorphic to the one associated to transposition. 
\item[(b)] Let ${}^{+}: M_{2h,2k}(\F) \rightarrow M_{2k,2h}(\F)$ be defined by:
 $a^{+}=H_{2k}a^tH_{2h}^{-1}$ for some skew-symmetric matrices 
$H_{2i}\in GL_{2i}(\F)$, 
then the corresponding $3$-bracket (\ref{antitranspose})
 is isomorphic to the one associated to map (\ref{symplectic}).  
\end{itemize}
\end{lemma}
{\bf Proof.} $(a)$ First of all notice that every symmetric matrix $A\in GL_n(\F)$ can be written as the product $B^{t}B$ where 
$B\in GL_n(\F)$
and $B^t$ is its transpose. Hence let $h=x^{t}x\in GL_m(\F)$ and $k=y^{t}y\in GL_n(\F)$ and consider the map $\varphi: M_{m,n}(\F) \rightarrow M_{m,n}(\F)$
defined by $\varphi(u)=x^{-1}uy$. For $a,b,c\in M_{m,n}(\F)$, let $[a,b,c]=ab^tc-cb^ta$ and $[a,b,c]^*=ab^*c-cb^*a$. Then we have:
$\varphi([a,b,c])=\varphi(ab^tc-cb^ta)=x^{-1}ab^tcy-x^{-1}cb^tay=
(x^{-1}ay)(y^{-1}b^tx)(x^{-1}cy)-(x^{-1}cy)(y^{-1}b^tx)(x^{-1}ay)=$
$\varphi(a)y^{-1}(y^{-1})^t(\varphi(b))^tx^tx\varphi(c)-\varphi(c)y^{-1}(y^{-1})^t(\varphi(b))^tx^tx\varphi(a)=
\varphi(a)k^{-1}(\varphi(b))^th\varphi(c)-
\varphi(c)k^{-1}(\varphi(b))^t$ $h\varphi(a)=[\varphi(a), \varphi(b), \varphi(c)]^*$,
and this shows that the $3$-brackets $[\,,\,,]$ and $[\,,\,,]^*$ are isomorphic.

$(b)$ Now let $H_{2k}=B^{-1}J_{2k}B$ and $H_{2h}=A^{-1}J_{2h}A$ for some 
matrices $A\in GL_{2k}$ and $B\in GL_{2h}$
such that $B^{-1}=B^t$, $A^{-1}=A^t$. 
For $a,b,c\in M_{2h,2k}(\F)$, let $[a,b,c]^*=ab^*c-cb^*a$ with
$b^*=J_{2k}b^tJ_{2h}^{-1}$. Let $\varphi: M_{2h,2k}(\F) \rightarrow M_{2h,2k}(\F)$ be defined
by $\varphi(x)=A^{-1}xB$. Then we have:
$\varphi([a,b,c]^*)=\varphi(ab^*c-cb^*a)=A^{-1}ab^*cB-A^{-1}cb^*aB=
(A^{-1}aB)(B^{-1}b^*A)(A^{-1}cB)-(A^{-1}cB)(B^{-1}b^*A)(A^{-1}aB)=$
$\varphi(a)B^{-1}J_{2k}b^tJ_{2h}^{-1}A\varphi(c)-\varphi(c)B^{-1}J_{2k}b^tJ_{2h}^{-1}A\varphi(a)=$
$\varphi(a)H_{2k}B^tb^tA$ $H_{2h}^{-1}\varphi(c)-\varphi(c)H_{2k}B^tb^tAH_{2h}^{-1}\varphi(a)=
$ $\varphi(a)(H_{2k}\varphi(b)^tH_{2h}^{-1})\varphi(c)-\varphi(c)(H_{2k}\varphi(b)^tH_{2h}^{-1})\varphi(a)=[\varphi(a), \varphi(b),$
$ \varphi(c)]^{+}$,
and this shows that the $3$-brackets $[\,,\,,]^*$ and $[\,,\,,]^{+}$ are 
isomorphic.
\hfill$\Box$

\medskip

\begin{example}\label{C(n)}\em
Let us consider the map $\psi: M_{1,2n}(\F)\rightarrow M_{2n,1}(\F)$, defined by:
$\psi(X ~Y)=(Y~ -X)^t$, for $X,Y\in M_{1,n}$. Then $M_{1,2n}$ with $3$-bracket
\begin{equation}
[A,B,C]=-AB^tC+CB^tA-C\psi(A)\psi(B)^t
\label{C(n)product}
\end{equation}
is an $N=6$ $3$-algebra, which we denote by $C^3(2n)$.
\end{example}

The N=6 $3$-algebras $A^3(m,n;t)$ and $C^3(2n)$ were introduced in \cite{BL2}. 

\begin{lemma}\label{osp(2,2n)}
Let $[\cdot, \cdot, \cdot]_*$ be the $3$-bracket on $M_{1,2n}$ defined by
$[a,b,c]_*=-a(kb^t\alpha^{-1})c+c(kb^t\alpha^{-1})a-c\psi(a)\alpha(\psi(b))^tk^{-1}$
 for some $\alpha\in\F^\times$
and some symmetric matrix 
$k\in Sp_{2n}(\F)$. Then $[\cdot, \cdot, \cdot]_*$ is
isomorphic to (\ref{C(n)product}).
\end{lemma}
{\bf Proof.} 
It
is convenient to identify $M_{1,2n}(\F)$ with the set of matrices of the form
$\left(\begin{array}{cc}
0 & 0\\
\ell & m
\end{array}\right)\in M_{2,2n}$, with $\ell,m\in M_{1,n}$. Under this 
identification,
for $Z\in M_{1,2n}(\F)$,
$\psi(Z)=J_{2n}Z^tJ_2^{-1}$.
Note that we can write 
$k=yy^t$ for some matrix $y\in Sp_{2n}(\F)$
\footnote{As E.\ Vinberg explained to us, if $\sigma$ is an anti-involution of
a connected reductive group $G$ and $S$ denotes its fixed point set in $G$,
then, by a well-known argument of Cartan, any element $k\in S$ can be
represented in the form $k=y\sigma(y)$ for some $y\in G$, provided that $S$
is connected.}, and, for $\left(\begin{array}{cc}
\alpha & 0\\
0 & \alpha^{-1}\end{array}\right)\in GL_2(\F)$, $h=xx^t$ for some 
$x\in GL_2(\F)$.
Consider the map $\varphi: M_{1,2n}(\F) \rightarrow M_{1,2n}(\F)$ defined by: 
$\varphi(u)=xuy^{-1}$
and let $[\cdot,\cdot,\cdot]$ be the $3$-bracket (\ref{C(n)product}). Then
we have:
$\varphi([a,b,c])=-(xay^{-1})yb^tx^{-1}(xcy^{-1})+(xcy^{-1})yb^tx^{-1}(xay^{-1})
-(xcy^{-1})y\psi(a)\psi(b)^ty^{-1}=$
$-\varphi(a)yb^tx^{-1}\varphi(c)+\varphi(c)yb^tx^{-1}\varphi(a)-
\varphi(c)y\psi(a)(\psi(b))^ty^{-1}
=-\varphi(a)k\varphi(b)^t$ $h^{-1}\varphi(c)+
\varphi(c)k\varphi(b)^th^{-1}\varphi(a)
-\varphi(c)\psi(\varphi(a))h\psi(\varphi(b))^tk^{-1}=
[\varphi(a),\varphi(b),$
$\varphi(c)]_*$,
since $y\psi(a)(\psi(b))^ty^{-1}$ $=\psi(\varphi(a))h\psi(\varphi(b))^tk^{-1}$.
Indeed we have: $\psi(\varphi(a))h\psi(\varphi(b))^tk^{-1}=
J_{2n}(xay^{-1})^tJ_2^{-1}xx^t(J_{2n}(xb$ $y^{-1})^tJ_2^{-1})^t(y^t)^{-1}y^{-1}$
$=J_{2n}(y^{-1})^ta^tx^tJ_2^{-1}xx^tJ_2xby^{-1}J_{2n}^{-1}(y^t)^{-1}y^{-1}=yJ_{2n}$
$a^tbJ_{2n}^{-1}y^{-1}
=y\psi(a)\psi(b)^t$ $y^{-1}$ since
$y^tJ_{2n}y=J_{2n}$ and $x^tJ_2^{-1}hJ_2x=I_2$.
\hfill$\Box$

\begin{example}\label{K(m,2)}\em
Let $P$ be a generalised Poisson algebra with bracket $\{.,.\}$ and derivation 
$D$ (see \cite{CantaK2} for the definition). 
Let $\sigma$ be a Lie algebra automorphism
of $P$ such that $-\sigma$ is an associative algebra automorphism, 
$\sigma^2=1$ and $\sigma\circ D=-D\circ \sigma$. Then
$P$ with the $3$-bracket:
\begin{equation}
[f,g,h]=\{f, \sigma(g)\}h+\{f,h\}\sigma(g)+f\{\sigma(g),h\}+
D(f)\sigma(g)h-f\sigma(g)D(h),
\label{eq:1.1}
\end{equation}
is an $N=6$ 3-algebra.

For example, consider the generalised Poisson algebra $P(m,0)$  in the (even) 
indeterminates $p_1, \dots,p_k, q_1, \dots, q_k$
(resp.\ $p_1, \dots,p_k, q_1, \dots, q_k, t$) if $m=2k$ (resp.\ $m=2k+1$)
endowed with the bracket: 
\begin{equation} 
\{f,g\}=(2-E)(f)\frac{\partial g}{\partial t}-\frac{\partial f}{\partial t}(2-E)(g)+\sum_{i=1}^k(\frac{\partial f}{\partial p_i}
\frac{\partial g}{\partial q_i}-
\frac{\partial f}{\partial q_i}\frac{\partial g}{\partial p_i}),
\label{Poisson}
\end{equation}
where $E=\sum_{i=1}^k(p_i\frac{\partial}{\partial p_i}+q_i\frac{\partial}
{\partial q_i})$
(the first two terms in (\ref{Poisson}) vanish if $m$ is even) and the 
derivation $D=2\frac{\partial}{\partial t}$ (which is 0 if $m$ is even).
Then the map 
$\sigma_{\varphi}: f(p_i,q_i)\mapsto -f(\varphi (p_i),\varphi(q_i))$ 
(resp. $\sigma_{\varphi}: f(t,p_i,q_i)\mapsto -f(\varphi(t),\varphi (p_i),\varphi(q_i))$, 
where $\varphi$ is an involutive linear change of variables 
(i.e. $\varphi^2=1$),
multiplying by $-1$ the 1-form  
$\sum_i (p_i dq_i - q_i dp_i)$ if $m$ is even 
(resp. $dt + \sum_i (p_i dq_i - q_i dp_i)$ if $m$ is odd), 
satisfies the conditions described above, hence
the corresponding $3$-bracket (\ref{eq:1.1})
defines on $P(m,0)$ an $N=6$ $3$-algebra structure. We denote this $3$-algebra
by $P^3(m;\varphi)$.
\end{example}

%

\begin{example}\label{S(1,2)new}\em
Let $A=\F[[x]]^{\langle 1\rangle}\oplus\F[[x]]^{\langle 2\rangle}$ be the 
direct sum of two copies of the algebra $\F[[x]]$; for $f\in\F[[x]]$, 
denote by $f^{\langle i\rangle}$ the
corresponding element in $\F[[x]]^{\langle i\rangle}$. Set $D=d/dx$ and let 
$a=(a_{ij})\in M_{2,2}(\F)$.
We define the following $3$-bracket on $A$ ($i,j=1$ or $2$):
$$[f^{\langle i\rangle},g^{\langle i\rangle},h^{\langle i\rangle}]=
(-1)^ia_{ij}((fD(h)-D(f)h)g(\varphi(x)))^{\langle i\rangle}
~~{\mbox{for}}~~j\neq i;$$
$$[f^{\langle i\rangle},g^{\langle j\rangle},h^{\langle i\rangle}]=
(-1)^ia_{jj}((fD(h)-D(f)h)g(\varphi(x)))^{\langle i\rangle}
~~{\mbox{for}}~~j\neq i;$$
$$[f^{\langle 1\rangle},g^{\langle j\rangle},h^{\langle 2\rangle}]=
a_{j1}((fD(g(\varphi(x)))-D(f)g(\varphi(x)))h)^{\langle 1\rangle})+
a_{j2}(f(hD(g(\varphi(x)))-D(h)g(\varphi(x))))^{\langle 2\rangle},$$
and extend it to $A$ by skew-symmetry in the first and third entries.
If $a\in SL_2(\F)$ and either $a^2=-1$ and $\varphi=-1$, or $a^2=1$ and 
$\varphi=1$, then $(A, [\cdot,\cdot,\cdot])$ is an 
$N=6$ $3$-algebra, which we denote by $SW^3(a)$. 
Note that if $a^2=1$, i.e., up to rescaling, $a=1$, and $\varphi=1$, then
we get the $N=8$ $3$-algebra $SW^3$ \cite{CantaK5}. 
\end{example}

\begin{example}\label{SKO(2,3;1)new}\em
Let $A=\F[[x_1,x_2]]$ and $D_i= \frac{\partial}{\partial x_i}$ for $i=1,2$.
We consider the following $3$-bracket:
\begin{equation}
[f,g,h]=\det
\left(
\begin{array}{ccc}
f & \varphi(g) & h\\
D_1(f) & D_1(\varphi(g)) & D_1(h)\\
D_2(f) & D_2(\varphi(g)) & D_2(h)
\end{array}\right)
\label{scircsign}
\end{equation}
where $\varphi$ is an automorphism of the algebra $A$.
If $\varphi$ is a linear change of variables with determinant equal to 1, and 
$\varphi^2=1$, then $A$ with $3$-bracket (\ref{scircsign}) is an $N=6$ 
$3$-algebra, which we denote by $S^3(\varphi)$. Note that if $\varphi=1$,
then we get the $N=8$ $3$-algebra $S^3$ \cite{CantaK5}. 
Clearly, all $3$-algebras $S^3(\varphi)$ with $\varphi \neq 1$ are isomorphic 
to each other.
\end{example}

\begin{example}\label{SHO(3,3)new}\em
Let $A=\F[[x_1,x_2,x_3]]$ and $D_i=\frac{\partial}{\partial x_i}$ for 
$i=1,2,3$.  Consider the following $3$-bracket on $A$:
\begin{equation}
[f,g,h]=\det
\left(
\begin{array}{ccc}
D_1(f) & D_1(\varphi(g)) & D_1(h)\\
D_2(f) & D_2(\varphi(g)) & D_2(h)\\
D_3(f) & D_3(\varphi(g)) & D_3(h)
\end{array}\right)
\label{scircfi}
\end{equation}
where $\varphi$ is an automorphism of the algebra $A$.
If $\varphi$ is a linear change of variables with determinant equal to 1, and 
$\varphi^2=1$, then $A$ with product (\ref{scircfi}) is an $N=6$ 
$3$-algebra, which we denote by $W^3(\varphi)$. Note that if $\varphi=1$,
then we get the $N=8$ $3$-algebra $W^3$ \cite{CantaK5}. Clearly, all 
$3$-algebras $W^3(\varphi)$ with $\varphi \neq 1$ are isomorphic to each other.
\end{example}

One can check directly that the above examples are indeed $N=6$ $3$-algebras.
However a proof of this without any computations will follow from the
connection of $N=6$ $3$-algebras to Lie superalgebras, discussed below.
The main result of the paper is the following theorem.
\begin{theorem}
The following is a complete list of simple linearly compact $N=6$
$3$-algebras over $\C$:

(a) finite-dimensional: 
$A^3(m,n;t)$, $A^3(2m,2n;st)$, $C^3(2n)$ ($m,n \geq 1$);

(b) infinite-dimensional: $P^3(m; \varphi)$ ($m \geq 1$),
$SW^3(a)$, $S^3(\varphi)$, $W^3(\varphi)$.
\end{theorem}

{\bf Proof.}
Theorem \ref{th:2.3} from Section 2 reduces the classification in question to that of 
the pairs $(L,\sigma)$, where
$L=L_{-1}\oplus L_0 \oplus L_1$
is a simple linearly compact Lie superalgebra with a consistent $\ZZ$-grading
and $\sigma$ is a graded conjugation of $L$. A complete list of possible such
$L=L_{-1}\oplus L_0 \oplus L_1$ is given by Remarks \ref{fdshortconsistent} and \ref{completelyodd} from Section \ref{section3}.
Finally, a complete list of graded conjugations of these $L$ is given by
Propositions \ref{gradedconjforfds} and \ref{listofgradedconj} from Section 4.

By Theorem \ref{th:2.3}$(b)$, the $N=6$ $3$-algebra is identified with $\Pi L_{-1}$,
on which the $3$-bracket is given by the formula 
$[a,b,c]=[[a,\sigma(b)],c]$. This formula, applied to the $\ZZ$-graded
Lie superalgebras $L$ with graded conjugations, described by Proposition \ref{gradedconjforfds}
$(a), (b), (c)$ in the finite-dimensional case produces the $3$-algebras
$A^3(m,n;t)$, $A^3(2m,2n;st)$, $C^3(2n)$, respectively, and those, described by
Proposition \ref{listofgradedconj}$(a)$ and $(b)$, $(c)$ and $(d)$, $(e)$ and $(f)$, $(g)$ 
in the infinite-dimensional 
case produces the $3$-algebras $P^3(m;\varphi)$, 
$SW^3(a)$, 
$S^3(\varphi)$, $W^3(\varphi)$,
respectively. The fact that all of them are indeed $N=6$ $3$-algebras follows
automatically from Theorem \ref{th:2.3}$(b)$. 
\hfill$\Box$

\section{Palmkvist's construction}
\begin{definition}
Let $\g=\oplus_{j\in\Z} \g_j$ be a 
Lie superalgebra with a consistent $\Z$-grading. A graded conjugation 
of $\g$ is a Lie
superalgebra automorphism
$\varphi: \g \rightarrow \g$ such that
\begin{enumerate}
\item $\varphi(\g_j)=\g_{-j}$
\item $\varphi^2(x)=(-1)^kx ~{\mbox{for}}~ x\in\g_k$.
\end{enumerate}
\end{definition}
\begin{theorem}
Let $\g=\oplus_{j\geq -1}\g_j$ be a $\Z$-graded consistent Lie superalgebra with a graded conjugation  $\varphi$. 
Then the $3$-bracket
$$[u,v,w]:=[[u,\varphi(v)],w]$$
defines on $\Pi\g_{-1}$ an $N=6$ 3-algebra structure.
\end{theorem}
{\bf Proof.} 
Since the grading of $\g$ is consistent, $\g_{-1}$ and $\g_1$ are completely odd and $\g_0$
is even. For $u,v,w\in \g_{-1}$ we thus have:
$[u,v,w]:=[[u,\varphi(v)],w]=[u,[\varphi(v),w]]=-[[\varphi(v),w],u]=-[[w,\varphi(v)],u]=-[w,v,u]$,
which proves property $(a)$ in Definition \ref{N6}.

Besides, for $u,v,x,y,z\in\g_{-1}$ we have:
$[u,v,[x,y,z]]-[x,y,[u,v,z]]=[[u,\varphi(v)],[[x,\varphi(y)],z]]-[[x,\varphi(y)],[[u,\varphi(v)],z]]=$
$[[u,\varphi(v)],[x,\varphi(y)]],z]+[[x,\varphi(y)],[[u,\varphi(v)],z]]-[[x,\varphi(y)],[[u,\varphi(v)],z]]=$
$[[u,\varphi(v)],[x,\varphi(y)]],z]=[[[u,\varphi(v)],x],\varphi(y)],z]+
[[[x,[[u,\varphi(v)],\varphi(y)]],z]=[[[u,\varphi(v)],x],\varphi(y)],z]-
[[[x,\varphi([[\varphi(u),v],y])],z]=$
$[[[u,v,x],y,z]-[x,[v,u,y],z]$
\hfill$\Box$

\bigskip

We shall now associate to an $N=6$ 3-algebra
$T$ with $3$-bracket $[\cdot, \cdot, \cdot]$, a $\Z$-graded Lie superalgebra 
$Lie T=Lie_{-1} T\oplus Lie_0 T\oplus Lie_1 T$, as follows. For $x,y\in T$, 
denote by $L_{x,y}$ the endomorphism
of $T$ defined by $L_{x,y}(z)=[x,y,z]$. 
Besides, for $x\in T$, denote by $\varphi_x$ the map in 
$Hom(\Pi T\otimes \Pi T, \Pi T)$
defined by $\varphi_x(y,z)=-[y,x,z]$. Here, as usual, $\Pi T$ denotes the 
vector space $T$ with odd parity.

We let $Lie_{-1} T=\Pi T$, $Lie_0 T=\langle L_{x,y} ~|~ x,y\in T\rangle$, 
$Lie_1 T=\langle \varphi_x ~|~ x\in T\rangle$,
and let $Lie T=Lie_{-1}T \oplus Lie_0 T \oplus Lie_1T$.  
Define the map $\sigma: Lie T \longrightarrow Lie T$ by 
$(x,y,z\in T)$:
$$z\mapsto -\varphi_z,
~~~\varphi_z\mapsto z, \,\,\,L_{x,y}\mapsto -L_{y,x}.$$
\begin{theorem}
\label{th:2.3}
(a) $Lie T$ is a $\Z$-graded Lie superalgebra with a short consistent grading,
satisfying the following two properties:

(i) any non-zero $\ZZ$-graded ideal of $Lie T$ has a non-zero 
intersection with both $Lie_{-1}T$ and $Lie_1 T$;

(ii) $[Lie_{-1} T, Lie_1 T]=Lie_0 T$.

(b) $\sigma$ is a graded conjugation of the $\ZZ$-graded Lie superalgebra
$Lie T$ and the 3-product on $T$ is recovered from the bracket on $Lie T$
by the formula $[x,y,z]=[[x,\sigma(y)],z]$.

(c) The correspondence $T\longrightarrow (Lie T, \sigma)$ 
is bijective and functorial between
$N=6$ $3$-algebras and the pairs $(Lie T,\sigma)$, where  
$Lie T$ is a $\Z$-graded Lie superalgebra with a short consistent grading,
satisfying properties (i) and (ii), and $\sigma$ is a graded conjugation
of $Lie T$.

(d) A $3$-algebra $T$ is simple (resp. finite-dimensional or linearly compact)
if and only if $Lie T$ is.  
\end{theorem}
{\bf Proof.} For $x,y,z\in T$, we have:
$[\varphi_x, z]=-L_{z,x}$, and
\begin{equation}
[L_{x,y},L_{x',y'}]=L_{[x,y,x'],y'}-L_{x',[y,x,y']}.
\label{eq:1.2}
\end{equation}
Note that
$[L_{x,y},\varphi_z]=-\varphi_{[y,x,z]}$.
Finally, $[[\varphi_x, \varphi_y], z]=[\varphi_x, [\varphi_y,z]]+[\varphi_y, [\varphi_x,z]]=
-[\varphi_x, L_{z,y}]-[\varphi_y, L_{z,x}]=[L_{z,y}, \varphi_x]+[L_{z,x}, \varphi_y]=
-\varphi_{[y,z,x]}-\varphi_{[x,z,y]}=0$. Hence $[\varphi_x, \varphi_y]=0$.

It follows that $Lie T$ is indeed a $\ZZ$-graded Lie superalgebra,
satisfying $(ii)$ and such that any non-zero ideal has a non-zero
intersection with $Lie_{-1}T$. It is straightforward to check $(b)$,
hence any non-zero ideal of $Lie T$ has a non-zero intersection
with $Lie_1T$, which completes the proof of $(a)$.

(c) is clear by construction. Since the simplicity of $T$, by definition,
means that all operators $L_{x,y}$ have no common non-trivial invariant 
subspace in $T$, it follows that $T$ is a simple $3$-algebra if and only if
$Lie_0 T$ acts irreducibly on $Lie_{-1} T$. Hence, by the properties $(i)$ and 
$(ii)$ of $Lie T$, $T$ is simple if and only if $Lie T$ is simple. The rest of 
(d) is clear as well. 
\hfill$\Box$

\begin{remark}\label{phi}\em 
$\sigma_{|Lie_0T}=1$ if an only if $T$ is a 3-Lie algebra. 
\end{remark}

\section{Short gradings and graded conjugations}\label{section3}
In this section we shall classify all short gradings of all simple linearly 
compact Lie superalgebras.
We shall describe the short gradings of the classical Lie superalgebras 
$\g$, $\g\neq Q(n), P(n)$, in terms of
linear functions $f$ on the set of roots of $\g$. For the description of the root systems of the classical Lie superalgebras, we shall refer to
\cite{K2}. 
As for the Lie superalgebra $Q(n)$, we will denote by
$e_i$ and $f_i$ the standard Chevalley generators of $Q(n)_{\0}$, and by $\bar{e}_i$ and $\bar{f}_i$ the
corresponding elements in $Q(n)_{\1}$.
Besides, we will identify $P(n)$ with the subalgebra of the Lie superalgebra 
$SHO(n,n)$ spanned by the following elements:
$\{x_ix_j, \xi_i\xi_j: i,j=1,\dots,n;~ x_i\xi_j: i\neq j=1,\dots,n;~
x_i\xi_{i}-x_{i+1}\xi_{i+1}: i=1, \dots,n-1\}$ (cf.\ \cite[\S 8]{CantaK4},
and thus describe the $\Z$-gradings of $P(n)$ as induced by the $\Z$-gradings of 
$SHO(n,n)$. Recall that the Lie superalgebra $W(0,n)$ is simple for $n\geq 2$, and for $n=2$ it is
isomorphic to the classical Lie superalgebra $osp(2,2)$ \cite{K1}. Moreover the Lie superalgebra $H(0,n)$ is simple
for $n\geq 4$, and for $n=4$ it is isomorphic to $psl(2,2)$.
Hence, when dealing with $W(0,n)$ (resp.\ $H(0,n)$) we shall always
assume $n\geq 3$ (resp.\ $n\geq 5$). Likewise, since 
$W(1,1)\cong K(1,2)$ and $S(2,1)\cong SKO(2,3;0)$ \cite[\S 0]{CantaK},
when dealing with $W(m,n)$ and $S(m,n)$ we shall always assume
$(m,n)\neq (1,1)$ and $(m,n)\neq (2,1)$, respectively.

\begin{proposition}\label{fdshort}
A complete list of simple finite-dimensional $\Z$-graded Lie superalgebras with a short grading
 $\g=\g_{-1}\oplus \g_0\oplus \g_1$
such that $\g_{-1}$ and $\g_1$ have the same dimension, 
is, up to isomorphism, as follows:
\begin{enumerate}
\item $A_m, \,B_m,\,C_m,\,D_m,\,E_6,\,E_7$ with the $\Z$-gradings, defined for each $s$ such that $a_s=1$ by
$f(\alpha_s)=1$, $f(\alpha_i)=0$ for all $i\neq s$, where $\alpha_1,\,\dots,
\alpha_m$ are simple roots and $\sum_i a_i\alpha_i$ is the highest root; 
\item $psl(m,n)$ with the $\Z$-gradings defined by:
$f(\epsilon_1)=\dots =f(\epsilon_k)=1$, 
$f(\epsilon_{k+1})=\dots =f(\epsilon_m)=0$,
$f(\delta_1)=\dots=f(\delta_h)=1$, 
$f(\delta_{h+1})=\dots=f(\delta_n)=0$,  for each $k=1,\dots, m$ and 
$h=1,\dots, n$;
\item $osp(2m+1,2n)$ with the $\Z$-grading defined by: $f(\epsilon_1)=1$, $f(\epsilon_i)=0$ for all $i\neq 1$,
$f(\delta_j)=0$ for all $j$;
\item $osp(2,2n)$ with the $\Z$-grading defined by:
$f(\epsilon_1)=1$, 
$f(\delta_j)=0$ for all $j$;
\item $osp(2,2n)$ with the $\Z$-grading defined by: $f(\epsilon_1)=1/2$, 
$f(\delta_j)=1/2$ for all $j$;
\item $osp(2m,2n)$, $m\geq 2$, with the $\Z$-grading defined by:
$f(\epsilon_1)=1$, $f(\epsilon_i)=0$ for all $i\neq 1$,
$f(\delta_j)=0$ for all $j$;
\item $osp(2m,2n)$, $m\geq 2$, with the $\Z$-grading defined by: $f(\epsilon_i)=1/2$, 
$f(\delta_j)=1/2$ for all $i, j$;
\item $D(2,1;\alpha)$ with the $\Z$-grading defined by: $f(\epsilon_1)=f(\epsilon_2)=1/2$, $f(\epsilon_3)=0$;
\item $F(4)$ with the $\Z$-grading defined by: $f(\epsilon_1)=1$, $f(\epsilon_i)=0$ for all $i\neq 1$, $f(\delta)=1$;
\item $Q(n)$ with the gradings defined by: $\deg(e_i)=\deg(\bar{e}_i)=-\deg(f_i)=-\deg(\bar{f}_i)=k_i$,
with $k_s=1$ for some $s$ and $k_i=0$ for all $i\neq s$;
\item $P(n)$ with $n=2h\geq 2$ and the gradings of type $(1,\dots,1,0,\dots,0|0,\dots,0,1,\dots,1)$ with $h$ 1's and $h$ 0's both
in the even and odd part; 
\item $H(0,n)$ with the grading of type $(|1,0,\dots,0,-1)$.
\end{enumerate}
\end{proposition}
{\bf Proof.} 
%
The claim for simple Lie algebras (the proof of which uses conjugacy of Borel 
subalgebras) is well known (see e.g. \cite{K2}).

In order to classify all short gradings of all classical Lie superalgebras 
$\g$,
we shall classify linear functions $f$ on the set of roots of $\g$ which take values $1$, $0$ or $-1$.

Let $\g=psl(m,n)$, i.e. $=sl(m,n)$ for $m\neq n$, and $=sl(n,n)/\F I_{n,n}$
for $m=n$. Since $\g$ has (even) roots $\pm(\epsilon_i-\epsilon_j)$ and $\pm(\delta_i-\delta_j)$,
 a linear function $f$ on the set of roots
taking values $0$ and $\pm 1$, is defined, up to a permutation of $\epsilon_i$'s and $\delta_j$'s, 
by: $f(\epsilon_1)=\dots =f(\epsilon_k)=a$, $f(\epsilon_{k+1})=\dots =f(\epsilon_m)=a-1$,
$f(\delta_1)=\dots=f(\delta_h)=b$, $f(\delta_{h+1})=\dots=f(\delta_n)=b-1$,  for some $a,b$, $k=0,\dots, m$, $h=0,\dots, n$.
On the other hand, $\g$ has (odd) roots $\pm(\epsilon_i-\delta_j)$, hence, either $b=a$ or $k=m$ and $b=a+1$. Since the value of $f$ on the roots is 
independent of $a$, we can let $a=1$.

Now let $\g=osp(2m+1,2n)$. Then $\g$ has roots $\delta_i$ and $2\delta_i$, hence $f(\delta_i)=0$ for every $i$. It follows that
for every $j$, either $f(\epsilon_j)=\pm 1$ or $f(\epsilon_j)=0$. If $f(\epsilon_j)=0$ for every $j$, then we get a grading which is
not short, hence we can assume, up to equivalence, that $f(\epsilon_1)=1$. Since $\g$ has roots $\epsilon_i\pm\epsilon_j$, it
follows that $f(\epsilon_j)=0$ for every $j\neq 1$.

Now let $\g=osp(2,2n)$. Then $\g$ has even roots of the form $\pm 2\delta_i$ and $\pm\delta_i\pm\delta_j$,
and odd roots of the form $\pm\epsilon_1\pm\delta_i$ hence, either $f(\delta_i)=0$ or $f(\delta_i)=\pm 1/2$.
If $f(\delta_i)=0$ for every $i$, then $f(\epsilon_1)=1$. If $f(\delta_k)=1/2$ for some $k$, then $f(\delta_i)=0$ for every $i$
and $f(\epsilon_1)=1/2$. We hence get two inequivalent short gradings.

Likewise, if $\g=osp(2m,2n)$ with $m\geq 2$, either $f(\delta_i)=0$ for every $i$, $f(\epsilon_k)=1$ for some $k$ and
$f(\epsilon_j)=0$ for every $j\neq k$, or $f(\delta_i)=0$ for every $i$ and $f(\epsilon_j)=1/2$ for every $j$.

Now let $\g=D(2,1;\alpha)$. Then $\g$ has roots $\pm 2\epsilon_i$, $i=1,2,3$, hence we may assume that $f(\epsilon_1)=1/2$. It
follows that, up to equivalence, $f(\epsilon_2)=1/2$ and $f(\epsilon_3)=0$.

Let $\g=F(4)$. Then $\g$ has even roots $\pm\epsilon_i\pm\epsilon_j$, $i\neq j$, $\pm\epsilon_i$, $\pm\delta$,  and
odd roots $1/2(\pm\epsilon_1\pm\epsilon_2\pm\epsilon_3\pm\delta)$. It follows that there exists some $k$ such that
$f(\epsilon_k)=1$, hence, $f(\epsilon_i)=0$ for every $i\neq k$ and $f(\delta)=\pm 1$, i.e., up to equivalence, $\g$ has
only one short grading.

Let $\g=G(3)$. Since $\delta$ and $2\delta$ are roots of $\g$, we have $f(\delta)=0$. Moreover, for every $i$,
 either $f(\epsilon_i)=0$ or $f(\epsilon_i)=\pm 1$. If $f(\epsilon_i)=0$ for every $i$, then the grading is not short, hence
 we may assume, up to equivalence, $f(\epsilon_1)=1$. It follows that, for $j=2,3$, either $f(\epsilon_j)=1$ or $f(\epsilon_j)=0$.
 But this contradicts the linearity of $f$ since $\epsilon_1+\epsilon_2+\epsilon_3=0$. Hence $G(3)$ has no short gradings.

Let $\g=Q(n)$. Then any $\Z$-grading of $\g$ is defined by setting 
$\deg(e_i)=\deg(\bar{e}_i)=-\deg(f_i)=-\deg(\bar{f}_i)=k_i\in\Z_{\geq 0}$.
Then, since $q(n)_{\bar{0}}\cong
A_n$, it is clear that such a grading has depth one if and only if all
$k_i$'s are 0 except for $k_s=1$ for some $s$.

Finally, we recall that a  complete list of $\Z$-gradings of depth 1 of all simple linearly compact Lie superalgebras
is given in \cite[Proposition 8.1]{CantaK4}. Then
if $\g$ is a simple finite-dimensional Lie superalgebra which is not classical, or $\g=P(n)$,
we select among these gradings the short ones such that $\g_1$ and $\g_{-1}$ have the same dimension, hence getting gradings
11. and 12. in the statement.
\hfill$\Box$

\begin{remark}\label{fdshortconsistent}\em It follows from Proposition 
\ref{fdshort} that
a complete list of simple finite-dimensional $\Z$-graded Lie superalgebras with a short consistent grading
 $\g=\g_{-1}\oplus \g_0\oplus \g_1$
such that $\g_{-1}$ and $\g_1$ have the same dimension, is, up
to isomorphism, as follows:
\begin{itemize}
\item[-] $psl(m,n)$ with $m,n\geq 1$, $m+n\geq 2$,
with the grading $f(\epsilon_1)=\dots=f(\epsilon_m)=1$,
$f(\delta_1)=\dots=f(\delta_n)=0$;
\item[-] $osp(2,2n)$, $n\geq 1$, with the grading $f(\delta_i)=0$ for all $i$,
$f(\epsilon_1)=1$.
\end{itemize}
Notice that $P(n)$ has, up to equivalence, a unique consistent $\Z$-grading, i.e., the
grading of type $(1,\dots,1|0,\dots,0)$. In this grading one has: 
$\g_{\0}\cong sl_{n+1}$, $\g_{-1}\cong \Lambda^2 (\F^{n+1})^*$ and $\g_1\cong S^2\F^{n+1}$, where $\F^{n+1}$ denotes
the standard $sl_{n+1}$-module, hence $\g_{1}$ and $\g_{-1}$ have different dimension.
Finally, no short grading of the Lie superalgebra $Q(n)$ is consistent since in this case $\g_{\1}$ is an
irreducible $\g_{\0}$-module.
\end{remark}

\begin{proposition}\label{gsshort}
A complete list of simple linearly compact infinite-dimensional $\Z$-graded 
Lie superalgebras with a short grading
 $\g=\g_{-1}\oplus \g_0\oplus \g_1$
such that $\g_{-1}$ and $\g_1$ have the same growth and the same size, is, up
to isomorphism, as follows:
\begin{itemize}
\item[-] $S(1,2)$ with the grading of type $(0|1,0)$;
\item[-] $S(1,2)$ with the grading of type $(0|1,1)$;
\item[-] $H(2k,n)$ with the grading of type $(0,\dots,0|1,0,\dots,0,-1)$;
\item[-] $K(2k+1,n)$ with the grading of type $(0,\dots,0|1,0,\dots,0,-1)$;
\item[-] $SHO(3,3)$ with the grading of type $(0,0,0|1,1,1)$;
\item[-] $SKO(2,3;\beta)$ with the grading of type $(0,0|1,1,1)$;
\item[-] $E(1,6)$ with the grading of type $(0|1,0,0,-1,0,0)$.
\end{itemize}
\end{proposition}
{\bf Proof.} All $\Z$-gradings of depth 1 of all infinite-dimensional linearly compact simple Lie superalgebras
are listed, up to isomorphism, in \cite[Proposition 9.1]{CantaK4}. Among these gradings we first select those
which are short, hence getting the following list:
\begin{enumerate}
\item[1)] $W(m,n)$, with $m>0$, $n\geq 1$, $(m,n)\neq (1,1)$, with the grading of type $(0,\dots,0|1,0,\dots,0)$;
\item[2)] $W(m,n)$, with $m>0$, $n\geq 1$, $(m,n)\neq (1,1)$, with the grading of type $(0,\dots,0|-1,0,\dots,0)$;
\item[3)] $S(m,n)$, with $m>1$ and $n\geq 1$, $(m,n)\neq (2,1)$, or $m=1$ and $n\geq 2$, with the grading of type $(0,\dots,0|1,0,\dots,0)$;
\item[4)] $S(1,2)$ with the grading of type $(0|1,1)$;
\item[5)] $H(2k,n)$ with the grading of type $(0,\dots,0|1,0,\dots,0,-1)$;
\item[6)] $K(2k+1,n)$ with the grading of type $(0,\dots,0|1,0,\dots,0,-1)$;
\item[7)] $SHO(3,3)$ with the grading of type $(0,0,0|1,1,1)$;
\item[8)] $SKO(2,3;\beta)$ with the grading of type $(0,0|1,1,1)$;
\item[9)] $E(1,6)$ with the grading of type $(0|1,0,0,-1,0,0)$.
\end{enumerate}
Let us consider $W(m,n)$ with $n\geq 1$ and the grading of type $(0,\dots,0|1,0,\dots,0)$. Then
$\g_{-1}=\langle \frac{\partial}{\partial\xi_1}\rangle\otimes\F[[x_1,\dots, x_m]]\otimes \Lambda(\xi_2,\dots,\xi_n)$
and 
$\g_{1}=\langle \xi_1\frac{\partial}{\partial x_i}, \xi_1\frac{\partial}{\partial \xi_j}
| i=1,\dots, m, j=2,\dots, n\rangle\otimes\F[[x_1,\dots,$ $x_m]]\otimes \Lambda(\xi_2,\dots,\xi_n)$. Therefore
$\g_{-1}$ and $\g_1$ have the same growth equal to $m$ but $\g_{-1}$ has size $2^{n-1}$ and
$\g_1$ has size $(m+n-1)2^{n-1}$. It follows that for $m>0$, $n\geq 1$ and $(m,n)\neq (1,1)$,
$\g_{-1}$ and $\g_1$ do not have the same size. Likewise case 2) is ruled out.

Now let us consider $S(m,n)$ with $n\geq 1$ and the grading of type $(0,\dots,0|1,0,\dots,0)$. Then
$\g_{-1}=\langle \frac{\partial}{\partial\xi_1}\rangle\otimes\F[[x_1,\dots, x_m]]\otimes \Lambda(\xi_2,\dots,\xi_n)$
and 
$\g_{1}=\{f\in \langle \xi_1\frac{\partial}{\partial x_i}, \xi_1\frac{\partial}{\partial \xi_j}
| i=1,\dots, m, j=2,\dots, n\rangle\otimes\F[[x_1,\dots, x_m]]\otimes \Lambda(\xi_2,\dots,\xi_n) ~|~ div(f)=0\}$. Therefore
$\g_{-1}$ and $\g_1$ have the same growth equal to $m$ but $\g_{-1}$ has size $2^{n-1}$ and
$\g_1$ has size $(m+n-2)2^{n-1}$. It follows that for $m>1$, $(m,n)\neq (2,1)$, or $m=1$ and  $n\geq 2$,
$\g_{-1}$ and $\g_1$ do not have the same size unless $m=1$ and $n=2$.
\hfill$\Box$

\begin{remark}\label{completelyodd}\em
It follows from Proposition 
\ref{gsshort} that
a complete list of simple infinite-dimensional $\Z$-graded Lie superalgebras with a short consistent grading
 $\g=\g_{-1}\oplus \g_0\oplus \g_1$
such that $\g_{-1}$ and $\g_1$ have the same growth and size, is, up
to isomorphism, as follows:
\begin{itemize}
\item[-] $S(1,2)$ with the grading of type $(0|1,1)$;
\item[-] $H(2k,2)$ with the grading of type $(0,\dots,0|1,-1)$;
\item[-] $K(2k+1,2)$ with the grading of type $(0,\dots,0|1,-1)$;
\item[-] $SHO(3,3)$ with the grading of type $(0,0,0|1,1,1)$;
\item[-] $SKO(2,3;\beta)$ with the grading of type $(0,0|1,1,1)$.
\end{itemize}
\end{remark}


\section{Classification of graded conjugations}
In this section we shall classify all graded conjugations $\sigma$ of all 
$\Z$-graded simple linearly compact Lie superalgebras $\g$ with a short 
consistent grading
$\g=\g_{-1}\oplus\g_0\oplus\g_1$. In Lemma \ref{generalform} and 
Proposition \ref{listofgradedconj} we shall assume that $\F=\C$.
This assumption can be removed with a little extra work.

\begin{remark}\label{1}\em
If $\varphi$ is an automorphism of $\g$ preserving the grading and $\sigma$ is a graded conjugation, 
then $\varphi\sigma\varphi^{-1}$
is again a graded conjugation which is equivalent to $\sigma$. Indeed, we have: 
$\varphi\sigma\varphi^{-1}\varphi\sigma\varphi^{-1}=\varphi\mu\varphi^{-1}=\mu$, where $\mu_{|\g_k}=(-1)^kid$.
\end{remark}

\begin{remark}\label{2}\em
If $\varphi$ is an involution of $\g$ preserving the grading and commuting with $\sigma$, then $\varphi\sigma$ is
again a graded conjugation.
\end{remark}

\begin{proposition}\label{gradedconjforfds}
The following is a complete list, up to equivalence, of graded conjugations of all simple finite-dimensional Lie superalgebras:
\begin{itemize}
\item[(a)] $\g=sl(m,n)/\F I\delta_{m,n}$: $\sigma_1\left(
\begin{array}{cc}
a & b\\
c & d
\end{array}\right)=\left(
\begin{array}{cc}
-a^t & c^t\\
-b^t & -d^t
\end{array}\right)$.
\item[(b)] $\g=sl(2h,2k)/\F I\delta_{2h,2k}$:
$\sigma_2\left(
\begin{array}{cc}
a & b\\
c & d
\end{array}\right)=\left(
\begin{array}{cc}
-a^{st} & c^{st}\\
-b^{st} & -d^{st}
\end{array}\right)$,
where $st$ denotes the symplectic transposition defined by (\ref{symplectic}).
\item[(c)] $\g=osp(2,2n)$: $\sigma_1$.
\end{itemize}
\end{proposition}
{\bf Proof.}
The Lie superalgebra $sl(m,n)$ has a short consistent grading such that 
$\g_0=\g_{\0}$ consists of
matrices of the form 
$\left(
\begin{array}{cc}
\alpha & 0\\
0 & \delta
\end{array}\right)$, where  
$tr \alpha=tr\delta$, $\g_{-1}$ is the set of matrices of the form
$\left(
\begin{array}{cc}
0 & 0\\
\gamma & 0
\end{array}\right)$, and
$\g_{1}$ is the set of matrices of the form
$\left(
\begin{array}{cc}
0 & \beta\\
0 & 0
\end{array}\right)$.

For $m\neq n$ every automorphism of $sl(m,n)$ is either of the form $Ad~diag(A, B)$ for some matrices $A\in GL_m(\F)$, $B\in GL_n(\F)$,
or of the form $Ad~diag(A, B)\circ \sigma_1$ \cite{S}. Note that
$\sigma_1^2=Ad~diag (I_m, -I_n)$ and $\sigma_1^3=\sigma_1^{-1}=Ad~diag (I_m,-I_n)\circ\sigma_1$. For $\left(
\begin{array}{cc}
a & b\\
c & d
\end{array}\right)\in sl(m,n)$, we have:
$Ad~diag(A,B)\left(
\begin{array}{cc}
a & b\\
c & d
\end{array}\right)=
\left(
\begin{array}{cc}
AaA^{-1} & AbB^{-1}\\
BcA^{-1} & BdB^{-1}
\end{array}\right)$, hence
 every automorphism $Ad~diag(A, B)$ maps $\g_{1}$ (resp.\ $\g_{-1}$) to itself and does not define
 a graded conjugation of $\g$.
Let $\varphi_{A,B}=Ad~diag(A,B)\circ \sigma_1$. Then
$\varphi_{A,B}\left(
\begin{array}{cc}
a & b\\
c & d
\end{array}\right)=
\left(
\begin{array}{cc}
-Aa^tA^{-1} & Ac^tB^{-1}\\
-Bb^tA^{-1} & -Bd^tB^{-1}
\end{array}\right)$.
If $\varphi_{A,B}$ is a graded conjugation of $\g$, then, by definition,
${\varphi_{A,B}}_{|\g_{\0}}^2=1$, hence $A^tA^{-1}=\lambda I_m$ and $B^tB^{-1}=\rho I_n$, for
some $\lambda, \rho\in \F$. Besides, since ${\varphi_{A,B}}_{|\g_{\1}}^2=-1$, we have 
$\rho=\lambda$.  It follows that $A^t=\lambda A$ hence, by transposing both sides of the equality,
$A=\lambda A^t=\lambda^2 A$, i.e., $\lambda^2=1$. Therefore,
either $A^t=A$ and $B^t=B$, or $A^t=-A$, $B^t=-B$ and $m$ and $n$ are even (since $A$ and $B$ are invertible).
The thesis then follows from Lemma \ref{sl(m,n)}.

If $m=n$, in addition to the automorphisms described above, $sl(n,n)/\F I_{n,n} $ has automorphisms of the form $Ad~diag(A, B)\circ \Pi$,
$Ad~diag(A, B)\circ \Pi\circ\sigma_1$ and $Ad~diag(A, B)\circ\sigma_1\circ \Pi$,
where $A, B\in GL_n(\F)$,
and $\Pi$ is defined as follows: for
$\left(
\begin{array}{cc}
a & b\\
c & d
\end{array}\right)\in sl(n,n)$, 
$\Pi \left(
\begin{array}{cc}
a & b\\
c & d
\end{array}\right)=\left(
\begin{array}{cc}
d & c\\
b & a
\end{array}\right)$  \cite{S}.
Note that $\sigma_1\circ\Pi\circ\sigma_1=\Pi$ and $\Pi\circ\sigma_1\circ\Pi=\sigma_1^{-1}$.
The automorphisms of the form $Ad~diag(A, B)\circ \Pi\circ\sigma_1$ and $Ad~diag(A, B)\circ\sigma_1\circ \Pi$
map $\g_{1}$ (resp.\ $\g_{-1}$) to itself, hence they do not define graded conjugations of $\g$.
Let $\psi_{A,B}=Ad~diag(A,B)\circ \Pi$. Then
$\psi_{A,B}\left(
\begin{array}{cc}
a & b\\
c & d
\end{array}\right)=
\left(
\begin{array}{cc}
AdA^{-1} & AcB^{-1}\\
BbA^{-1} & BaB^{-1}
\end{array}\right)$. It follows that ${\psi_{A,B}}_{|_{\g_{\0}}}^2=1$ if and only if $AB=BA=\lambda I_n$. As a consequence,
${\psi_{A,B}}_{|_{\g_{\1}}}^2=1$, hence $\psi_{A,B}$ does not define a graded conjugation.

The Lie superalgebra $osp(2,2n)$ has a short consistent grading such that $\g_{\0}$ consists of matrices of the form
$\left(
\begin{array}{cc}
A & 0\\
0 & D
\end{array}\right)$, where
$A=\left(
\begin{array}{cc}
\alpha & 0\\
0 & -\alpha
\end{array}
\right)$, $\alpha\in\F$, and $D$ lies in the Lie algebra $sp(2n)$, defined by
$J_{2n}$,
 $\g_{-1}$ is the set of matrices of the form
$\left(
\begin{array}{c|c}
0 & {\begin{array}{cc}
0 & 0\\
a & b
\end{array}}\\
\hline
{\begin{array}{cc}
b^t & 0\\
-a^t & 0
\end{array}} & 0
\end{array}\right)$, and
$\g_{1}$ is the set of matrices of the form
$\left(
\begin{array}{c|c}
0 & {\begin{array}{cc}
a & b\\
0 & 0
\end{array}}\\
\hline
{\begin{array}{cc}
0 & b^t\\
0 & -a^t 
\end{array}} & 0
\end{array}\right)$, with $a,b\in M_{1,n}$.

Every automorphism of $osp(2,2n)$ is either of the form $Ad~diag(A, B)$ for some matrices $A=diag(\alpha,\alpha^{-1})$,
$\alpha\in\F^{\times}$, $B\in Sp_{2n}(\F)$,
or of the form $Ad~diag(A, B)\circ \sigma_1$ \cite{S}. One can easily check that every automorphism of the form
$Ad~diag(A, B)$ sends $\g_{-1}$ (resp.\ $\g_1$) to itself, hence it does not define a graded conjugation of $\g$.
Let $\Phi_{A,B}=Ad~diag(A, B)\circ\sigma_1$. Then, using the same arguments as for the automorphisms $\varphi_{A,B}$
of the Lie superalgebra $sl(m,n)$, one can show that $\Phi_{A,B}$ defines a graded conjugation of $osp(2,2n)$ if
and only if $B$ is a symmetric matrix. Then the result follows from Lemma \ref{osp(2,2n)}. 
\hfill$\Box$

\begin{remark}\label{automorphisms}\em
If $\g$ is a simple infinite-dimensional linearly compact Lie superalgebra, then $Aut~\g$ contains a maximal reductive subgroup
which is explicitely described in \cite[Theorem 4.2]{CK2}. We shall denote this subgroup by $G$.
We point out that any reductive subgroup
of $Aut~\g$ is conjugate into $G$, in particular any finite order element of
$Aut~\g$ is conjugate to an element of $G$.
\end{remark}

\begin{example}\label{K(1,2)}\em
The grading of type $(0,\dots,0|1,-1)$ of $\g=H(2k,2)$ (resp.\ $K(2k+1,2)$) is short.
Let $A=\F[[p_1,\dots,p_k,q_1,\dots,q_k]]$ (resp.\ $A=\F[[t,p_1,\dots,p_k,q_1,\dots,q_k]]$). We have:

\medskip

$\g_{-1}= \langle \xi_2\rangle\otimes A$,

$\g_0=(\langle 1,\xi_1\xi_2\rangle\otimes A)/\F 1$ (resp.\ $\langle 1,\xi_1\xi_2\rangle\otimes A$),

$\g_{1}=\langle \xi_1\rangle\otimes A$.

\medskip

\noindent
For every linear involutive change of variables $\varphi$,
multiplying by $-1$ the 1-form
$\sum_{i=1}^k (p_i dq_i - q_i dp_i)$
(resp. $dt + \sum_{i=1}^k (p_i dq_i - q_i dp_i)$),
the following map is a graded conjugation of $\g$:
\begin{equation}
\begin{array}{ll}
f(p_i, q_i)\mapsto -f (\varphi(p_i), \varphi(q_i)) & (\mbox{resp.}~ f(t,p_i, q_i)\mapsto -f(\varphi(t), \varphi(p_i), \varphi(q_i))\\
f(p_i, q_i)\xi_1\xi_2\mapsto -f(\varphi(p_i), \varphi(q_i))\xi_1\xi_2 & (\mbox{resp.}~ f(t,p_i, q_i)\xi_1\xi_2\mapsto -f(\varphi(t), \varphi(p_i), \varphi(q_i))\xi_1\xi_2)\\
f(p_i, q_i)\xi_1\mapsto f (\varphi(p_i), \varphi(q_i))\xi_2 & (\mbox{resp.}~ f(t,p_i, q_i)\xi_1\mapsto f(\varphi(t), \varphi(p_i), \varphi(q_i))\xi_2)\\
f(p_i, q_i)\xi_2\mapsto -f (\varphi(p_i), \varphi(q_i))\xi_1 & (\mbox{resp.}~ f(t,p_i, q_i)\xi_2\mapsto -f(\varphi(t), \varphi(p_i), \varphi(q_i))\xi_1).
\end{array}
\label{gradedforH}
\end{equation}
\end{example}


\begin{example}\label{S(1,2)}\em
Let $\g=S(1,2)$, $SHO(3,3)$, or $SKO(2,3;1)$. Then the algebra of outer derivations of $\g$ contains $sl_2=
\langle e,h,f\rangle$, with $e=\xi_1\xi_2\frac{\partial}{\partial x}$ and $h=\xi_1\frac{\partial}{\partial\xi_1}+\xi_2\frac{\partial}{\partial\xi_2}$ if $\g=S(1,2)$, $e=\xi_1\xi_3\frac{\partial}{\partial x_2}-\xi_2\xi_3\frac{\partial}{\partial x_1}-\xi_1\xi_2\frac{\partial}{\partial x_3}$
and $h=\sum_{i=1}^3\xi_i\frac{\partial}{\partial\xi_i}$ if
$\g=SHO(3,3)$, $e=\xi_1\xi_2\tau$ and $h=1/2(\tau-x_1\xi_1-x_2\xi_2)$ if $\g=SKO(2,3;1)$. Let us denote by $G_{out}$
the subgroup of $Aut~\g$ generated by $\exp(ad(e))$, $\exp(ad(f))$ and $\exp(ad(h))$. We recall that 
$G_{out}\subset G$, where $G$ is the subgroup of
$Aut~\g$ introduced in Remark \ref{automorphisms} \cite[Remark 2.2, Theorem 4.2]{CK2}.
We shall denote by $U_-$ the one parameter group of automorphisms $\exp(ad(tf))$,  and by
$G_{inn}$ the subgroup of $G$ consisting of inner automorphisms. 
Finally, 
$H$ will denote the subgroup of
$Aut~\g$  consisting of  invertible changes of variables multiplying the volume form
(resp.\ the even supersymplectic form) by a constant if $\g=S(1,2)$ (resp.\ $\g=SHO(3,3)$),
or the odd supercontact form by a function if $\g=SKO(2,3;1)$ (see \cite[Theorem 4.5]{CK2}). 

The gradings of type $(0|1,1)$, $(0,0,0|1,1,1)$ and $(0,0|1,1,1)$ of $\g=S(1,2)$,
$SHO(3,3)$ and $SKO(2,3;1)$, respectively, are short, and the subspaces $\g_i$'s are as follows:

\medskip

\noindent
$\g=S(1,2)$:

$\g_{-1}=\langle \frac{\partial}{\partial\xi_1}, \frac{\partial}{\partial\xi_2}\rangle\otimes\F[[x]]$

$\g_0=\{f\in \langle \frac{\partial}{\partial x}, \xi_i\frac{\partial}{\partial\xi_j} | i,j=1,2\rangle\otimes\F[[x]], div(f)=0\}$

$\g_{1}=\{f\in \langle \xi_i\frac{\partial}{\partial x}, \xi_1\xi_2\frac{\partial}{\partial\xi_i} | i=1,2\rangle\otimes\F[[x]], div(f)=0\}$.

\medskip

\noindent
$\g=SHO(3,3)$:

$\g_{-1}=\F[[x_1,x_2,x_3]]/\F 1$

$\g_0=\{f\in \langle\xi_1, \xi_2, \xi_3\rangle\otimes\F[[x_1,x_2,x_3]] | \Delta(f)=0\}$

$\g_1=\{f\in \langle\xi_i\xi_j, i,j=1,2,3\rangle\otimes\F[[x_1,x_2,x_3]] | \Delta(f)=0\}$.

\medskip

\noindent
$\g=SKO(2,3;1)$:

$\g_{-1}=\F[[x_1,x_2]]$

$\g_0=\{f\in \langle\xi_1, \xi_2, \tau\rangle\otimes\F[[x_1,x_2]] | div_1(f)=0\}$

$\g_1=\{f\in \langle\tau\xi_i, \xi_1\xi_2 ~|~ i=1,2\rangle\otimes\F[[x_1,x_2]] | div_1(f)=0\}$.

\medskip

\noindent
In all these cases the map $s=\exp(ad(e))\exp(ad(-f))\exp(ad(e))$ is a graded conjugation of $\g$: for 
$z\in\g_{-1}$, $s(z)=[e,z]$; for $z\in\g_{1}$, $s(z)=-[f,z]$, for $z\in \g_0$, $s(z)=z$.
Note that each of the above gradings can be extended to $Der ~\g=\g\rtimes {\mathfrak a}$, with ${\mathfrak a}\supset sl_2$,
so that $e$ has degree 2, $h$
has degree 0, and $f$ has degree $-2$.
\end{example}

%
%
%
%
%
%

\begin{lemma}\label{generalform} 
If $\g$ is one of the following $\Z$-graded Lie superalgebras:
\begin{enumerate}
\item $S(1,2)$ with the grading of type $(0|1,1)$,
\item $SKO(2,3;1)$ with the grading of type
$(0,0|1,1,1)$, 
\item $SHO(3,3)$ with the grading of type $(0,0,0|1,1,1)$,
\end{enumerate} 
and $\sigma$ is a graded conjugation of $\g$, then $\sigma$
is conjugate to an automorphism of the form $s\circ\exp(ad(th))\circ\varphi$,
for some $t\in\F$ and some $\varphi\in G_{inn}$  such that $\varphi^2=1$.
\end{lemma}
{\bf Proof.} Let us first assume $\g=S(1,2)$ with the grading of type $(0|1,1)$, or
$\g=SKO(2,3;1)$ with the grading of type $(0,0|1,1,1)$. By \cite[Remark 4.6]{CK2}, if $\psi$ is an automorphism of $\g$ lying in $G$, then either $\psi\in U_-H\cap G$ or $\psi\in U_-sH\cap G$.
Note that $U_{-}H\cap G=U_-(H\cap G)$ and $U_{-}sH\cap G=U_-s(H\cap G)$, since $U_{-}\subset G$ and $s\in G$ \cite[Theorem 4.2]{CK2}. 
Here $H\cap G$ is the subgroup of $Aut~\g$ generated by $\exp(ad(e))$, $\exp(ad(h))$ and $G_{inn}$. Note that $G_{inn}\subset \exp(ad(\g_0))$.
Let $\psi\in U_-(H\cap G)$.  
Then $\psi=\exp(ad(tf))\psi_0$ for some $t\in\F$ and some $\psi_0\in H\cap G$. For $x\in\g_1$, we have:
$\psi(x)=\exp(ad(tf))(\psi_0(x))=\psi_0(x)+t[f,\psi_0(x)]$, since $\psi_0(x)\in \g_1$. In particular, $\psi(x)\notin\g_{-1}$.
Now let $\sigma$ be a graded conjugation of $\g$. Then we may assume, 
up to conjugation, that $\sigma$ lies in $G$. Since $\sigma$
exchanges $\g_1$ and $\g_{-1}$, by the observation above 
$\sigma\in U_{-}s(H\cap G)$, i.e., $\sigma=\exp(ad(tf))s\varphi_0\varphi_1$
for some $t\in\F$, some $\varphi_1\in G_{inn}$ and some $\varphi_0$ lying in the subgroup generated by
$\exp(ad(h))$ and $\exp(ad(e))$. 
We can assume $\varphi_0=\exp(ad(\beta e))\exp(ad(\alpha h))$ for some $\alpha, \beta\in\F$, i.e., 
$\sigma=\exp(ad(tf))s\exp(ad(\beta e))\exp(ad(\alpha h))\varphi_1$. Since $\varphi_1(\g_1)=\g_1$ and
$\varphi_1(\g_{-1})=\g_{-1}$, for $x\in \g_1$ we have:
$$\sigma(x)=-\exp({\alpha})[f, \varphi_1(x)]=s\circ\exp(ad(\alpha h))\circ\varphi_1(x).$$
For $x\in\g_{-1}$ we have:
\begin{equation}
\sigma(x)=\exp({-\alpha})([e, \varphi_1(x)]+(t-\beta)\varphi_1(x)).
\label{sigma}
\end{equation}
Notice that if $x\in\g_{-1}$, then $[e, \varphi_1(x)]\in\g_1$ and $\varphi_1(x)\in\g_{-1}$. Since
$\sigma(\g_{-1})=\g_1$, we have 
$$(\ref{sigma})=\exp({-\alpha})[e, \varphi_1(x)]=s\circ\exp(ad(\alpha h))\circ\varphi_1(x).$$
Therefore $\sigma=s\circ\exp(ad(\alpha h))\circ\varphi_1$.
Now notice that $s\circ\exp(ad(\alpha h))=\exp(ad(-\alpha h))\circ s$, $\exp(ad(\alpha h))\circ\varphi_1=\varphi_1\circ\exp(ad(\alpha h))$,
and $s\circ\varphi_1=\varphi_1\circ s$. It follows that
$\sigma^2=s\circ\exp(ad(\alpha h))\circ\varphi_1\circ s\circ\exp(ad(\alpha h))\circ\varphi_1=
s\circ\exp(ad(\alpha h))\exp(ad(-\alpha h))\circ s\circ\varphi_1^2=s^2\circ\varphi_1^2$, therefore
$\varphi_1^2=1$.

Now let $\g=SHO(3,3)$ with the grading of type $(0,0,0|1,1,1)$. As in the previous cases, if $\sigma$ is a graded conjugation of $\g$,
then $\sigma\in U_-s(H\cap G)$. Here $H\cap G$ is the subgroup of $Aut~\g$ generated by $\exp(ad(e))$, $\exp(ad(h))$, $\exp(ad(\Phi))$, and $G_{inn}$,
where $\Phi=\sum_{i=1}^3(-x_i\frac{\partial}{\partial x_i}+\xi_i\frac{\partial}{\partial \xi_i}$) and
$G_{inn}$ is generated by $\exp(ad(x_i\xi_j))$ with $i, j=1,2,3$, $i\neq j$, and is thus isomorphic to $SL_3$. Note that 
$\Phi$ commutes with $G_{inn}$ and $h$. We may hence assume $\sigma=\exp(ad(tf))s\varphi_0\varphi_1$, for some $t\in\F$, some
$\varphi_1\in G_{inn}$ and some $\varphi_0$ lying in the subgroup generated by $\exp(ad(e))$, $\exp(ad(h))$, $\exp(ad(\Phi))$,
i.e., $\varphi_0=\exp(ad(\gamma e))\circ\exp(ad(\beta \Phi))\circ\exp(ad(\alpha h))$, for some $\alpha, \beta, \gamma\in\F$. Arguing as for $S(1,2)$ and $SKO(2,3;1)$,
one shows that, in fact, $\sigma=s\circ\exp(ad(\beta \Phi))\circ\exp(ad(\alpha h))\circ\varphi_1$.
Besides, the following commutation relations hold: $s\circ\exp(ad(\alpha h))=\exp(ad(-\alpha h))\circ s$,
$\exp(ad(\alpha h))\circ\varphi_1=\varphi_1\circ\exp(ad(\alpha h))$, $s\circ\varphi_1=\varphi_1\circ s$,
$\exp(ad(\beta\Phi))\circ s_{|\g_{-1}}=\exp({3\beta})s\circ\exp(ad(\beta\Phi))_{|\g_{-1}}$,
$\exp(ad(\beta\Phi))\circ s_{|\g_{1}}=\exp({-3\beta})s\circ\exp(ad(\beta\Phi))_{|\g_{1}}$.
It follows that, since $\sigma^2=s^2$,we have:
$\exp({3\beta})\exp(ad(\beta\Phi))^2{\varphi_1^2}_{|\g_{-1}}=1$,
$\exp({-3\beta})\exp(ad(\beta\Phi))^2{\varphi_1^2}_{|\g_{1}}=1$.
Note that $G_{inn}$ acts on $\g_{-1}$ by the standard action of vector fields on functions. In particular $V=\langle x_1, x_2, x_3\rangle$
is stabilized by this action and $\exp(ad(\Phi))$ acts on $V$ by scalar multiplication by $\exp({-1})$. It follows that if $F\in SL_3$ is the
matrix of the action of $\varphi_1$ on $V$, then $\exp({\beta})I_3F^2=I_3$, hence $\exp({3\beta})=\exp({3\beta})\det(F)^2=1$. It follows
that $(\exp(ad(\beta\Phi))\circ \varphi_1)^2=1$, moreover, we can assume $\beta=0$, since $\exp(\beta) I_3\in SL_3$.
Hence  $\sigma=s\circ \exp(ad(\alpha h))\circ \varphi_1$, for some $\alpha\in\F$ and some $\varphi_1\in G_{inn}$ such that
$\varphi_1^2=1$.
\hfill$\Box$

\begin{proposition}\label{listofgradedconj}
The following is a complete list, up to equivalence, of graded conjugations of 
all simple infinite-dimensional linearly
compact Lie superalgebras $\g$:
\begin{itemize}
\item[a)] $\g=H(2k,2)$: $\sigma$ is the automorphism of $\g$ defined by (\ref{gradedforH}).
\item[b)] $\g=K(2k+1,2)$: $\sigma$ is the automorphism of $\g$ defined by (\ref{gradedforH}).
\item[c)] $\g=S(1,2)$: $s$.
\item[d)] $\g=S(1,2)$: $\sigma=s\circ\exp(ad(\alpha h'))\circ\varphi_0$, where $\exp({2\alpha})=-1$, 
$h'=2x\frac{\partial}{\partial x}+\xi_1\frac{\partial}{\partial \xi_1}+\xi_2\frac{\partial}{\partial \xi_2}$, 
$\varphi_0$ lies in the $SL_2$-subgroup of $G_{inn}$, and $\varphi_0^2=-1$.
\item[e)] $\g=SKO(2,3;1)$: $s$.
\item[f)] $\g=SKO(2,3;1)$: $\sigma=s\circ\exp(ad(\alpha h'))$, where $\exp({2\alpha})=-1$ and $h'=x_1\xi_1+x_2\xi_2+\tau$.
\item[g)] $\g=SHO(3,3)$: $\sigma=s\circ\varphi$ with $\varphi\in G_{inn}$ such that $\varphi^2=1$.
\end{itemize}
\end{proposition}
{\bf Proof.} By definition of graded conjugation, $\g$ is, up to 
isomorphism, one of the $\Z$-graded Lie superalgebras listed in Remark
\ref{completelyodd}.
 Let $\g=H(2k,2)$ with the grading of type $(0,\dots,0|1,-1)$ (see Example \ref{K(1,2)}), and let $\sigma$
be a graded conjugation of $\g$. By Remark \ref{automorphisms} we can assume that
$\sigma\in G=\F^{\times}(Sp_{2k}\times O_2)$ \cite[Theorem 4.2]{CK2}. Note
that $G$ consists of linear changes of variables preserving the symplectic form
up to multiplication by a non-zero scalar. Since $\sigma$ exchanges $\g_{-1}$ and $\g_1$, and ${\sigma^2}_{|\g_{-1}}=-1$,
we have: $\sigma(\xi_1)=a\xi_2$ and $\sigma(\xi_2)=-\frac{1}{a}\xi_1$, for some $a\in\F^\times$, hence, up to multiplication by a scalar,
we may assume that $a=1$. It follows that, for $f\in\F[[p_i,q_i]]$,
$\sigma(f)=\sigma([\xi_1, f\xi_2])=[\xi_2, \sigma(f\xi_2)]=[\xi_2, \tilde{f}\xi_1]=\tilde{f}$ for some
$\tilde{f}\in\F[[p_i,q_i]]$, i.e., $\sigma(f\xi_2)=\sigma(f)\xi_1$. Likewise,
$\sigma(f\xi_1)=-\sigma(f)\xi_2$, and $\sigma(f\xi_1\xi_2)=\sigma(f)\xi_1\xi_2$. Besides,
for $f,g\in\F[[p_i,q_i]]$
$\sigma([f\xi_1, g\xi_2])=\sigma(fg+[f,g]\xi_1\xi_2)=-[\sigma(f)\xi_2, \sigma(g)\xi_1]=
-\sigma(f)\sigma(g)+[\sigma(f), \sigma(g)]\xi_1\xi_2$. Hence $\sigma(fg)=-\sigma(f)\sigma(g)$, i.e.,
$-\sigma$ is an automorphism of $\F[[p_i, q_i]]$ as an associative algebra.
It follows that $\sigma$ is defined as follows:
\begin{equation}
\begin{array}{l}
f(p_i, q_i)\mapsto -f (\varphi(p_i), \varphi(q_i))\\
f(p_i, q_i)\xi_1\xi_2\mapsto -f (\varphi(p_i), \varphi(q_i))\xi_1\xi_2\\
f(p_i, q_i)\xi_1\mapsto f (\varphi(p_i), \varphi(q_i))\xi_2\\
f(p_i, q_i)\xi_2\mapsto -f (\varphi(p_i), \varphi(q_i))\xi_1
\end{array}
\end{equation}
for some linear change of even variables $\varphi$. Since $\sigma(\xi_1)=\xi_2$
and $\sigma(\xi_2)=-\xi_1$,
$\sigma$ multiplies the odd part $d\xi_1d\xi_2$ of the symplectic form by $-1$,
hence $\varphi$ multiplies the even part of the symplectic form by $-1$.
Moreover, $\sigma^2(f(p_i, q_i))=f(\varphi^2(p_i), \varphi^2(q_i))$, hence $\varphi^2=1$.
This concludes the proof of a).
The same arguments prove b). In this case, one has $\sigma\circ\frac{\partial}{\partial t}=-\frac{\partial}{\partial t}\circ \sigma$.

Let $\g=SKO(2,3;\beta)$ with $\beta\neq 1$ and the grading of type $(0,0|1,1,1)$.  By \cite[Theorem 4.2]{CK2},
$G$ is generated by $\exp(ad(\tau+x_1\xi_1+x_2\xi_2))$ and $G_{inn}$. Note that $G_{inn}$ is contained 
in $\exp(ad(\g_0))$, hence no automorphism of $\g$ exchanges $\g_1$ and $\g_{-1}$. It follows that $\g$ has no
graded conjugations. 

Let $\g=S(1,2)$ and let $\sigma$ be a graded conjugation of $\g$. Then, by Lemma \ref{generalform}, 
$\sigma=s\circ\exp(ad(th))\circ\varphi$ for some $t\in\F$ and some $\varphi\in G_{inn}$ such that
$\varphi^2=1$. The group $G_{inn}$ is generated by $\exp(ad(h'))$, $\exp(ad(\xi_1\frac{\partial}{\partial\xi_2}))$,
$\exp(ad(\xi_2\frac{\partial}{\partial\xi_1}))$ and $\exp(ad(\xi_1\frac{\partial}{\partial\xi_1}-\xi_2\frac{\partial}{\partial\xi_2}))$,
hence we may write $\varphi=\exp(ad(\alpha h'))\varphi_0$ for some $\alpha\in\F$ and some $\varphi_0$ in the $SL_2$-subgroup of $G_{inn}$
generated by $\exp(ad(\xi_1\frac{\partial}{\partial\xi_2}))$,
$\exp(ad(\xi_2\frac{\partial}{\partial\xi_1}))$ and $\exp(ad(\xi_1\frac{\partial}{\partial\xi_1}-\xi_2\frac{\partial}{\partial\xi_2}))$.
Note that $\varphi(\frac{\partial}{\partial x})=\exp({-2\alpha})\frac{\partial}{\partial x}$, therefore $\exp({2\alpha})=\pm 1$ since $\varphi^2=1$.
Besides, if $z\in\g_{-1}=\langle \frac{\partial}{\partial \xi_1}, \frac{\partial}{\partial \xi_2}\rangle\otimes\F[[x]]$, then
$\varphi^2(z)=\exp(-2\alpha)\varphi_0^2(z)$, therefore 
either 

$i)$ $\exp({2\alpha})=1$ and ${\varphi_0}_{|\g_{-1}}^2=1$;

\noindent
or

$ii)$ $\exp({2\alpha})=-1$ and ${\varphi_0}_{|\g_{-1}}^2=-1$.

\noindent
In case $i)$ we have ${\varphi_0}_{|\g_{-1}}=\pm 1$, since $\varphi_0\in SL_2$. Then ${\varphi_0}_{|\g_{0}}=1$ since $\varphi_0$ acts
on $\g_0$ via the adjoint action. It follows that $\sigma_{|\g_{0}}=1$, since  $\sigma=s\circ\exp(ad(th))\circ\exp(ad(\alpha h'))\varphi_0$,
$\exp(ad(\alpha h'))_{|\g_{0}}=1$ since $\exp({2\alpha})=1$, $\exp(ad(th))_{|\g_{0}}=1$ and
$s_{|\g_{0}}=1$. By Remark \ref{phi} and the classification of $3$-Lie algebras
obtained in \cite{CantaK5}, we conclude that $\sigma$ is conjugate to $s$.

In case $ii)$ $\varphi_0$ corresponds to a $2\times 2$
matrix  of the form
$\left(\begin{array}{cc}
a & b\\
c & -a
\end{array}\right)$ such that $a^2+bc=-1$. The corresponding element $\sigma=s\circ\exp(ad(th))\circ\exp(ad(\alpha h'))\varphi_0$
acts on $\g_{-1}$ as follows:
$$x^r\frac{\partial}{\partial\xi_1}\mapsto \exp({2r\alpha})\exp({-t-\alpha})(-arx^{r-1}\xi_1\xi_2\frac{\partial}{\partial\xi_1}-
brx^{r-1}\xi_1\xi_2\frac{\partial}{\partial\xi_2}+ax^r\xi_2\frac{\partial}{\partial x}-bx^r\xi_1\frac{\partial}{\partial x});$$
$$x^r\frac{\partial}{\partial\xi_2}\mapsto \exp({2r\alpha})\exp({-t-\alpha})(-crx^{r-1}\xi_1\xi_2\frac{\partial}{\partial\xi_1}+
arx^{r-1}\xi_1\xi_2\frac{\partial}{\partial\xi_2}+cx^r\xi_2\frac{\partial}{\partial x}+ax^r\xi_1\frac{\partial}{\partial x}).$$
It follows that, up to rescaling, we may assume that $t=0$ hence getting $d)$.

In order to classify the graded conjugations of $\g=SKO(2,3;1)$ we argue in a similar way as for $S(1,2)$. Namely, let $\sigma$ be
such a map, then, by Lemma \ref{generalform}, $\sigma=s\circ \exp(ad(th))\circ\varphi$ for some $t\in\F$ and some
$\varphi\in G_{inn}$ such that $\varphi^2=1$. The group $G_{inn}$ is generated by $\exp(ad(h'))$, $\exp(ad(x_1\xi_2))$, $\exp(ad(x_2\xi_1))$
and $\exp(ad(x_1\xi_1-x_2\xi_2))$, hence we may write $\varphi=\exp(ad(\alpha h'))\varphi_0$ for some $\alpha\in\F$ and
some $\varphi_0$ in the $SL_2$-subgroup of $G_{inn}$ generated by  $\exp(ad(x_1\xi_2))$, $\exp(ad(x_2\xi_1))$
and $\exp(ad(x_1\xi_1-x_2\xi_2))$.
Note that this $SL_2$-subgroup
acts on $\g_{-1}=\F[[x_1, x_2]]$ via the standard action of vector fields on functions, and stabilizes the subspaces of
$\F[[x_1, x_2]]$ consisting of homogeneous polynomials of fixed degree.
We have: $\varphi(1)=\exp(ad(\alpha h'))(1)=\exp({-2\alpha})$, therefore, since
$\varphi^2=1$, either $\exp({2\alpha})=1$
or $\exp({2\alpha})=-1$. If $\exp({2\alpha})=1$, then,
 for $f\in\g_{-1}$, $\varphi(f)=\varphi_0(f)$; if $\exp({2\alpha})=-1$, then, for $f\in\g_{-1}$, $\varphi(f)=-\varphi_0(f^-)$, where $f^-(x_1, x_2)=f(-x_1, -x_2)$. It follows that
${\varphi_0^2}_{|\g_{-1}}=1$. 
In particular,
if $V=\langle x_1, x_2\rangle$, then ${\varphi_0}_{|V}=\pm 1$, i.e., $\varphi_0=\exp(ad(A(x_1\xi_1-x_2\xi_2))$ with $A\in\F$ such that
$\exp(A)=\pm 1$. It follows that $\sigma=s\circ\exp(ad(th))\circ\exp(ad(\alpha h'))\exp(ad(A(x_1\xi_1-x_2\xi_2))$,
for some $t, \alpha, A\in\F$ such that $\exp({2\alpha})=\pm 1$, $\exp(A)=\pm 1$. We now consider the restriction of $\sigma$
to $\g_0$. We have: $\sigma_{|\g_0}=\exp(ad(\alpha h'))\exp(ad(A(x_1\xi_1-x_2\xi_2))_{|\g_0}$, since
$s_{|\g_0}=1$ and $\exp(ad(th))_{|\g_0}=1$. It is then easy to check that
either 

\medskip

$i)$ $\exp({2\alpha})\exp(A)=1$ and $\sigma_{|\g_0}=1$;

\noindent
or 

\medskip

$ii)$ $\exp({2\alpha})\exp(A)=-1$.

\medskip
\noindent
In case $i)$, by Remark \ref{phi} and the classification of $3$-Lie algebras
obtained in \cite{CantaK5}, $\sigma$ is conjugate to $s$.
In case $ii)$, $\sigma=s\circ\exp(ad(th))\circ\exp(ad(\alpha h'))\exp(ad(A(x_1\xi_1-x_2\xi_2))$ acts on $\g_{-1}=\F[[x_1, x_2]]$ as follows:
$$f\mapsto -\exp({-t})s(f^-), ~~~\mbox{if}~~ \exp(A)=1, \exp({-2\alpha})=-1;$$
$$f\mapsto \exp({-t})s(f^-), ~~~\mbox{if}~~ \exp(A)=-1, \exp({-2\alpha})=1.$$
Therefore, changing the sign if necessary, we may assume that we are in the first case and in this case we may assume, up to rescaling, that
$A=0=t$, hence getting $f)$.

Let $\g=SHO(3,3)$ with the grading of type $(0,0,0|1,1,1)$. Then, by Lemma \ref{generalform}, $\sigma=s\circ\exp(ad(th))\circ\varphi$
for some $\varphi\in G_{inn}$ such that $\varphi^2=1$, and some $t\in\F$. For $a\in\g_{-1}$, we have:
$$\sigma(a)=s(\exp(ad(th))(\varphi(a))=\exp({-t})s(\varphi(a))=\exp({-t})[e,\varphi(a)],$$
since $\varphi(a)\in \g_{-1}$.
Up to rescaling, we can thus assume that $t=0$, hence getting $g)$. 
Note that in this case $G_{inn}\cong SL_3$.
\hfill$\Box$

%
%
%
%
%
%
%
%
%
%
%
%
%
\begin{remark}\em
It is proved in \cite{CantaK5} that there are no simple linearly compact 
$N=8$ $3$-superalgebras, which are not $3$-algebras. 
On the contrary, there are many simple linearly compact $N=6$ $3$-
superalgebras beyond $3$-algebras. We are planning to classify them
in a subsequent publication.  
\end{remark}

\section{$N=5$ $3$-algebras}

Based on the discussion in \cite{BB}, the following seems to be a right 
definition of an $N=5$ 3-algebra.

\begin{definition}\label{N5} An $N=5$ 3-algebra is a 3-algebra whose $3$-bracket  $[\cdot,\cdot,\cdot]$ satisfies
the following axioms:
\begin{enumerate}
\item[$(a)$] $[u,v,w]=[v,u,w]$
        \item[$(b)$]  $[u,v,[x,y,z]]=[[u,v,x],y,z]+[x,[u,v,y],z]+[x,y,[u,v,z]]$
        \item[$(c)$] $[u,v,w]+[v,w,u]+[w,u,v]=0$.
\end{enumerate}
\end{definition}

The following example is inspired by \cite{J}: we just replace 
$\ZZ /2 \ZZ$-graded Lie algebras by Lie superalgebras.

\begin{example}\em
Let $\g$ be a Lie superalgebra. Define, for
$a,b,c \in V:=\Pi \g_{\1}$, $[a,b,c]=[[a,b],c]$. Then 
$V$ with this $3$-bracket is an $N=5$ 3-algebra.
Indeed, $(a)$ follows from the skew-commutativity of the (super)bracket, and 
$(b)$ and $(c)$ from the
(super) Jacobi identity.

Conversely, any $N=5$ $3$-algebra can be constructed in this way. Namely, let 
$(N, [\cdot,\cdot,\cdot])$
be an $N=5$ $3$-algebra. Set $\g_{\1}(N)=\Pi N$ and let $\g_{\0}(N)$ be the subalgebra
of the Lie superalgebra $End(N)$ spanned by elements $L_{a,b}$, with $a,b\in N$, defined by:
$$L_{a,b}(c)=[a,b,c].$$
Note that, by property $(a)$ of Definition \ref{N5}, $L_{a,b}=L_{b,a}$. Moreover, $[L_{a,b}, L_{c,d}]=L_{[a,b,c],d}+L_{[a,b,d],c}$. Let $\g(N)=\g_{\0}(N)+\g_{\1}(N)$
with $[L_{a,b},c]=L_{a,b}(c)=-[c, L_{a,b}]$, and $[a,b]=L_{a,b}$, for $a,b,c\in\g_{\1}(N)$.
Then $\g(N)$ is a Lie superalgebra.
Indeed, the skew-commutativity of the bracket follows immediately from the construction. Besides, the (super) Jacobi identity for $\g(N)$ can be proved
as follows: for $a,b,c,d,x\in\g_{\1}(N)$, 
$$\begin{array}{ll}
[[L_{a,b}, L_{c,d}],x] &=L_{[a,b,c],d}(x)+L_{[a,b,d],c}(x)=[[a,b,c],d,x]+[[a,b,d],c,x]\\
 & =
[a,b,[c,d,x]]-[c,d,[a,b,x]]=[L_{a,b}, [L_{c,d},x]]-[L_{c,d}, [L_{a,b},x]],
 \end{array}$$
where we used property $(b)$ of Definition \ref{N5}; besides,
$$[[L_{a,b}, c], d]=[[a,b,c],d]=L_{[a,b,c],d}=[L_{a,b}, L_{c,d}]-L_{c,[a,b,d]}=
[L_{a,b}, [c,d]]-[c, [L_{a,b}, d]].$$
Finally, $[[a,b],c]=[L_{a,b}, c]=[a,b,c]=
-[b,c,a]-[a,c,b]=[a,[b,c]]+[b,[a,c]]$, by property $(c)$ of Definition \ref{N5}.
\end{example}

A skew-symmetric bilinear form $(.,.)$ on a finite-dimensional
$N=5$ $3$-algebra is called invariant
if the $4$-linear form $([a,b,c],d)$ on it is invariant
under permutations $(ab)$, $(cd)$ and $(ac)(bd)$ (which generate a dihedral group 
of order 8).

It is easy to see that if $\g$ is a finite-dimensional Lie superalgebra,
then the restriction of any invariant supersymmetric bilinear form 
$(.,.)$ on $\g$ to $\g_{\bar{1}}$ 
defines on the $N=5$ $3$-algebra 
$V=\Pi \g_{\bar{1}}$ an invariant bilinear form.
If, in addition, $\g$ is a simple Lie superalgebra and the
bilinear form is non-degenerate, then $\g$ is isomorphic to one of
Lie superalgebras $psl(m,n)$, $osp(m,n)$, $D(2,1;\alpha)$, $F(4)$,
$G(3)$, or $H(2k)$ \cite {K1}. All examples of the corresponding $N=5$
3-algebras appear in \cite{BB}, except for $\g=H(2k)$. In the latter case the
corresponding $N=5$ $3$-algebra is the subspace of odd elements
of the Grassmann algebra in $2k$ indeterminates $\xi_i$ with reversed parity,
endowed with the following $3$-bracket: $[a,b,c]=\{\{a,b\},c\}$, where
$\{a,b\}=\sum_i\frac{\partial a}{\partial \xi_i}\frac{\partial b}
{\partial \xi_i}$,
the invariant bilinear form being $(a,b)=$ coefficient of 
$\xi_1 ... \xi_{2k}$ in $ab^*$, where $b^*$ is the Hodge dual of $b$.
 
We show, in conclusion, how to associate an $N=5$ 3-algebra to an 
$N=6$ 3-algebra.
Let $(L,[\cdot,\cdot,\cdot]_6)$ be an $N=6$ 3-algebra. Let $T=L+L'$, 
where $L'=\langle\varphi_x, ~x\in L\rangle$,
$\varphi_x(y,z)=-[y,x,z]_6$. Let $\sigma : T \rightarrow T$ be
defined by:
$\sigma(z)=-\varphi_z, \sigma(\varphi_z)=z$ (cf.\ Remark \ref{phi}). 
Then $\sigma^2=-1$.
Now define on $T$ the following $3$-bracket ($a,b,c \in T$):
$$[a,b,c]_5=0 ~{\mbox{if}}~ a,b\in L ~{\mbox{or}}~ a,b\in L';$$
$$[a,b,c]_5=[b,a,c]_5=[a,\sigma(b),c]_6=[[a,b],c] ~{\mbox{if}}~ a\in L, b\in L'.$$
Then $(T,[\cdot,\cdot,\cdot]_5)$ is an $N=5$ $3$-algebra.


\end{document}